
\baselineskip=14pt
\parskip=10pt
\def\Tilde{\char126\relax}

\font\eightrm=cmr8  
\font\eighttt=cmtt8
\magnification=\magstephalf

\parindent=0pt
\overfullrule=0in
 
\bf
\centerline
{
The COMBINATORIAL ASTROLOGY of Rabbi ABRAHAM IBN EZRA
}
\rm
\bigskip
\centerline{ {\it 
Doron ZEILBERGER
}\footnote{$^1$}
{\eightrm  \raggedright
Department of Mathematics, Temple University,
Philadelphia, PA 19122, USA. 
{\eighttt zeilberg@math.temple.edu} \hfill \break
{\eighttt http://www.math.temple.edu/\Tilde zeilberg/   .}
Aug. 25, 1998. Supported in part by the NSF.
Written by Invitation of Jose L. Fernandez, who
promised to translate this into Spanish and publish
it in {\eighttt La Gaceta}, the Spanish analog of the Monthly, of which he
is editor. The English version is exclusive for
the Personal Journal of Ekhad and Zeilberger,
{\eighttt http://www.math.temple.edu/\Tilde zeilberg/pj.html} .
} 
}
 
Rabbi Abraham Ben Meir Ibn Ezra (1089-1164), from Tudela,
Spain, was one of the all-time greatest in the
following categories of  intellectual endeavor.
 
1) {\bf Poet}: e.g.: {\it aha yarad/al sfrad/ra min hashamaim/
eini eini yorda mayim} [Aha befell, on Spain, evil from the
heavens, my eye, my eye, water falling]).
 
2) {\bf Philosopher}: he was an extreme pantheist and neo-Platonist,
who influenced Spinoza in his abstract conception of God.
 
3) {\bf Biblical exegist}: he introduced a critical and
grammatical approach to biblical commentary, and was much
more `modern' than apologists like Rashi.
I particularly like his
commentary on Leviticus 18:20, where he said that
the sexual act can be divided into three kinds:
the good kind, for procreation, the bad one, to satisfy
beastly lust, and a neutral one, `to alleviate the
weariness of the body'.
 
4) {\bf Grammarian}: by writing in Hebrew,
he revealed to European scholars, who could
not read the Arabic literature, the ideas of his
predecessors, including the concept of the
three-letter {\it shoresh} (root), and the method
of {\it pealim} (verbs, actions).
 
5) {\bf Mathematician}: e.g. he popularized Zero 
(that he called {\it galgal}(wheel)), and wrote a very
influential textbook, `The Book of Number', where
among many other things, he showed how to 
reduce multiplication to taking squares
by using $(a+b)(a-b)=a^2-b^2$, and how to take squares
by using the recursions:
$(3a)^2=10a^2-a^2, (n \pm 1)^2= (n)^2 \pm 2n +1$,
 
6) {\bf Puzzler}: he (allegedly) saved his
life and the life of his disciples,
by solving the generalized Josephus problem
(there were 15 good guys and 15 bad guys on a boat,
when a storm started to rage,
15 passengers had to be thrown to the sea; how to
arrange the 30 people in a circle, by drowning every
ninth man, in such a way that the scoundrels
all drown?).

But what he was probably most proud of
was his {\bf ASTROLOGY}. And indeed,
he was the greatest astrologer of his time,
and probably in the `top three' of all times
(the two others being Ptolemy and Kepler).
In fact, a large part of his mathematics was inspired by
applications to astrology. 
 
One of his astrological treatises is called
{\it Sefer HaOlam} (The Book of the World). 
Written in a polemical style, its main message is to warn users
against 
`wrong' applications of astrology. Of course, like
most scholars until `modern' times, he was an ardent
believer in astrology, but only when it is practiced correctly.
 
In particular he warned that all the astronomical tables predicting
the times of planetary conjunctions are erroneous, because they
assume uniform motion of the planets. He also made a very good
point on the accumulation of errors, and for the need to account for
experimental errors, and how unreasonable it is to extrapolate from
ancient data. Hence, he only relied on astronomical observations
made by contemporary `sages of experiments'.
 
Ibn Ezra also knew how to compute the seventh row of
what later became to be called Pascal's triangle.
Except for the trivial $k=0,1$, he showed how to compute ${{7} \choose {k}}$ 
for $k\leq 7$.
 
The practical problem that inspired these calculations was to find
the number of possible planetary conjunctions.
As every educated person knows, there are exactly seven planets:
Sun, Moon, Mercury, Venus, Mars, Jupiter, and Saturn.
Whenever a subset of cardinality larger than one shows up at
the same sign, this has great astrological significance.
How many such events are possible?
 
Here is a translation, from Hebrew, of the first half of the 
relevant passage
in Ibn Ezra's Sefer HaOlam. He called planets `servants'
({\it meshartim}), which perhaps meant servants of God.
 
\it
And the combinations are hundred and twenty. And thus you can
know their number. [It is] Known that every calculation that adds
from one to any number that one wills, you can obtain by its value
[multiplied] by its half together with half of one, and here is
an example, we wanted to know what is the sum of the numbers from
one to twenty. We will multiply twenty by its half and half of
one, and [we get] two hundreds and ten. And now we can start
to know how many combinations [involving] two servants. And it is
known that the number of servants is seven. And Saturn can
combine with six other servants. And six by its half and
the half of one is one and twenty. And thus is the number of
combinations of twos. [Now] we wanted to know the number of
combinations of threes. Here we put Saturn and Jupiter 
and one of the others, their number is five. We multiply
five by two and a half and a half, and get fifteen \dots.
\rm
 
Ibn Ezra, then computes the number of conjunctions of
three planets without Saturn, and repeats the same argument
to get ${{5} \choose {2}}$, then, in turn, 
${{4} \choose {2}}$, ${{3} \choose {2}}$, and 
${{2} \choose {2}}$, totaling $35$.

What Ibn Ezra is doing here is using the formulas (in our notation)
$$
{{n} \choose {2}}= \sum_{i=1}^{n-1} i = 
(n-1) \left ( {{n-1} \over {2}}+{{1} \over {2}} \right )
\quad, 
\eqno(1)
$$
and
$$
{{n} \choose {k}}=\sum_{m=k-1}^{n-1}  {{m} \choose {k-1}}
\eqno(2)
$$.
 
For $k=3,4$, he uses $(2)$ repeatedly until he can use
$(1)$.
Translating his reasoning for the calculation of
${{7} \choose {4}}$ to our notation reads as follows.
$$
{{7} \choose {4}}=
{{6} \choose {3}}+{{5} \choose {3}}+{{4} \choose {3}}+
{{3} \choose {3}}=
$$
$$
[{{5} \choose {2}}+{{4} \choose {2}}+
{{3} \choose {2}}+{{2} \choose {2}}]+
[{{4} \choose {2}}+{{3} \choose {2}}+
{{2} \choose {2}}]+
[{{3} \choose {2}}+{{2} \choose {2}}]+
{{2} \choose {2}}]=
$$
$$
[ 4(4/2+1/2)+3(3/2+1/2)+2(2/2+1/2)+1(1/2+1/2)]+
[ 3(3/2+1/2)+2(2/2+1/2)+1(1/2+1/2)]+
$$
$$
[2(2/2+1/2)+1(1/2+1/2)]+1(1/2+1/2)=35 \quad .
$$
 
For $k=5,6,7$, Ibn Ezra only uses $(2)$, without $(1)$, invoking 
direct enumeration.

{\bf Morals} 
 
1) We are lucky to have modern notation.
Ibn Ezra was obviously far smarter than any of us,
yet he had to struggle so much because he lacked the
right notation and the precise idea of induction.
It was Rabbi Levi Ben Gerson, in 1321, who {\it rigorously proved}
the explicit expressions for the binomial
coefficients, and even he, almost two hundred years later,
had to go through a verbal
nightmare because he did not quite have modern algebraic
notation. I believe that we are about to witness
another notational revolution, inspired by programming
languages, in which future proofs would be phrased.
I am sure that our grandchildren will view our
current style of writing mathematics and proofs using
English or Spanish with
the same bewilderment and mild amusement that we view
Ibn Ezra and Levi Ben Gerson's mouthfuls.
 
2) Do not be superstitious about
so-called superstitions. Not only Abraham Ibn Ezra and
Levi Ben Gerson, but also Kepler and Newton, considered
Astrology as a real science, not a pseudo-science.
Who knows which of our current `scientific' views
will be considered superstition and hogwash by future
generations? I have two candidates: the actual
infinity, and the insistence on `rigorous proof'.

\bye